\documentclass[letterpaper, 10 pt, conference]{ieeeconf}
\usepackage{amsmath,amssymb,latexsym,float,epsfig,algorithm,subcaption,algorithmicx,algpseudocode,tabularx,pdflscape,epigraph,xparse,color,cite,url}
\usepackage{chngcntr,color}

\IEEEoverridecommandlockouts
\begin{document}
\title{\LARGE \bf On the Benefits of Surrogate Lagrangians\\ in Optimal Control and Planning Algorithms}

\author{Gerardo~De~La~Torre~and~Todd~D.~Murphey 
\thanks{G. De La Torre and T.~D. Murphey are with the Department of Mechanical Engineering, McCormick School of Engineering and Applied Science, Northwestern University, 2145 Sheridan Road, Evanston, IL 60208, USA. 
E-mail: gerardo.delatorre@northwestern.edu, t-murphey@northwestern.edu.}
\thanks{This work was supported by Army Research Office grant W911NF-14-1-0461.}
}

\maketitle

\begin{abstract}
This paper explores the relationship between numerical integrators and optimal control algorithms.
Specifically, the performance of the differential dynamical programming (DDP)  algorithm is examined when a variational integrator and a newly proposed surrogate variational integrator are used to propagate and linearize system dynamics. 
Surrogate variational integrators, derived from backward error analysis, achieve higher levels of accuracy while maintaining the same integration complexity as nominal variational integrators.
The increase in the integration accuracy is shown to have a large effect on the performance of the DDP algorithm.    
In particular, significantly more optimized inputs are computed when the surrogate variational integrator is utilized. 
\end{abstract}
\section{Introduction}
Due to the intractable nature of the Hamilton-Jacobi-Bellman equation numerical solvers for the nonlinear optimal control problem have been widely developed and implemented. 
Generally, numerical optimal control solvers are categorized as direct or indirect methods\cite{rao2009survey}. 
While indirect methods use the calculus of variations to obtain multiple-point boundary value problems, direct methods discretize the optimal control problem.
How the optimal control problem is discretize is an important aspect of direct methods.
In this paper, the differential dynamical programming (DDP) algorithm is used as a motivating example to explore the importance of discretization in direct methods. 
In particular, this paper examines the role of the selected numerical integrator used to propagate and linearize system dynamics in the DDP algorithm.

The DDP algorithm generates optimal control policies utilizing a quadratic approximation of the cost-to-go function and linearized state space dynamics\cite{DDPbook,DDPconstaint001,DDPminimax,SDDP}. 
Iterative linear quadratic regulators (iLQR) were derived similarly\cite{ilqr001}. 
Further analysis and development of the algorithm has addressed state and control constraints and proven, under mild assumptions, quadratic convergence of the solution\cite{DDPconstaint001, DDPconstaint002, DDPminimax,ShoemakerConv}.  
The algorithm has been successfully implemented in simulation and in a real-time aerospace application\cite{DLT:DISS}. 
Furthermore, the benefits of using variational integrators in DDP algorithms have been demonstrated \cite{IMAROB}.

Variational integrators are ideal candidates to used in an optimal control solver due to their long-term energy preserving properties and available structured linearization\cite{ober2011discrete,VI:Marsdem:01,VI:Murphey:02,hussein2006discrete,johnson2012trajectory}. 
However, traditional variational integrators can only achieve second-order accuracy and, therefore, their achievable benefits are limited.
This paper introduces the concept of surrogate Lagrangians for variational integrators. 
Surrogate variational integrators, derived from backward error analysis, achieve higher levels of accuracy while maintaining the same integration scheme as a nominal variational integrator \cite{reich1999backward,hairer1994backward,WARMING1974159,griffiths1986scope}.
Therefore, the proposed integrator maintains the geometric and structural properties of nominal variational integrators.

This paper explores the relationship between numerical integrators and optimal control algorithms.
Specifically, the performance of the DDP algorithm is examined when a variational integrator and a newly proposed surrogate variational integrator are used to propagate and linearize system dynamics. 
Numerical experiments demonstrate significantly more optimized inputs are computed when the surrogate variational integrator is utilized.
Therefore, the accuracy of numerical integrators may have a large impact on all algorithms used for control and state estimation. 

The organization of this paper is as follows. 
Section II gives an overview of variational integrators, structured linearization, and the differential dynamic programming algorithm.
The concepts of backward error analysis and surrogate Lagrangians are outlined in Section III.
The benefits of surrogate Lagrangians are demonstrated through a series of numerical experiments in Section IV.
Conclusions and future research directions are discussed in Section~V. 
\section{Background}
\subsection{Overview of Variational Integrators}
Variational integrators are used to propagate a dynamical system's (potentially forced and constrained) configuration through time. 
Rather than using a discrete approximation of the Euler-Lagrange equation, variational integrators are created by approximating the Lagrange-d'Alembert principle\cite{West:thesis,VI:Marsdem:01,VI:Murphey:01}.
The resulting approximation gives an implicit two-step mapping from two consecutive discrete configurations to the next $(q_{k-1},q_k)\rightarrow(q_{k+1})$.
Therefore,  the continuous system configuration is approximated by a sequence of discrete configuration vectors $(q_1,q_2,\dots,q_n)$.
Formally developing a variational integrator begins with defining the \emph{discrete Lagrangian} as
\begin{align}
L_{\textrm{d}}(q_k,q_{k+1})&= hL(\frac{q_{k+1}+q_k}{2},\frac{q_{k+1}-q_{k}}{h}),\nonumber\\
&\approx\int^{t_{k+1}}_{t_k}L(q(\tau),\dot{q}(\tau)) \ \textrm{d}\tau, \label{dis:Lag}
\end{align}
where $L(q(t),\dot{q}(t))$ is the Lagrangian of the dynamical system defined as
\begin{align}
L(q,\dot{q}) = T(q(t),\dot{q}(t)) - V(q(t)),
\end{align}
such that $T(q,\dot{q})$ and $V(q)$ are the system's kinetic and potential energy, respectively, and 
$q(t)$ is the system's configuration. 
Furthermore, $q_k$ is the system's configuration at time $t_k$ and $h$ is the discretization time step ($h = t_{k+1}-t_{k}$). 
The definition of the discrete Lagrangian is not unique and a generalized midpoint approximation can be used to derive other discrete Lagrangians  such that $L_{\textrm{d}}(q_k,q_{k+1}) = hL(q_{k+1}(1-\alpha)+q_k\alpha,\frac{q_{k+1}-q_{k}}{h})$ where $0\leq\alpha\leq1$.
However, by setting $\alpha=0.5$ second-order convergence of the integrator is obtained\cite{West:thesis}. 
Likewise, \emph{left and right discrete forces} are defined as
\begin{align}
&F^{-}_{\textrm{d}}(q_k,q_{k+1},u_{k})\hspace{-0.05cm}\cdot\hspace{-0.05cm} \delta q_{k} \hspace{-0.1cm}+ \hspace{-0.1cm}F^{+}_{\textrm{d}}(q_k,q_{k+1},u_{k})\hspace{-0.05cm}\cdot\hspace{-0.05cm} \delta q_{k+1}\nonumber\\
&\quad\approx\int^{t_{k+1}}_{t_k}F(q(\tau),\dot{q}(\tau),u(\tau))\cdot \delta q \ \textrm{d}\tau,
\end{align}
where 
\begin{align}
&F^{\pm}_{\textrm{d}}(q_k,q_{k+1},u_{k})= \frac{h}{2}F(\frac{q_{k+1}+q_k}{2},\frac{q_{k+1}-q_{k}}{h}, u_k),
\end{align}
$u(t)$ is the system's input (if any), and $\delta q$ is a virtual displacement. 
The defined quantities are used to obtain a discrete approximation of the \emph{Lagrange-d'Alembert} principle \cite{2002analytical}\footnote{The slot derivative $D_jG(y_1,y_2,\dots)$ is defined to be the derivative of $G$ with respect to its ${j^\textrm{th}}$ argument. That is, $D_jG(y_1,y_2,\dots) = \frac{\partial G}{\partial y_j}$.}:
\begin{align}
\delta S[q(t)] &= \delta \int^{t_\textrm{f}}_{0} L(q,\dot q) \ \textrm{d}\tau + \int^{t_\textrm{f}}_{0}F(q,\dot{q},u)\cdot{\delta q} \ \textrm{d}\tau,\nonumber\\
&\approx \sum^{n-1}_{k=1}\big(D_1L_{\textrm{d}}(q_k,q_{k+1}) + D_2L_{\textrm{d}}(q_{k-1},q_{k})\big)\cdot \delta q_{k} \nonumber\\
&\quad+ \sum^{n-1}_{k=1}\big(F^{-}_{\textrm{d}}(q_k,q_{k+1},u_k)\nonumber\\ &\quad\quad+ F^{+}_{\textrm{d}}(q_{k-1},q_{k},u_{k})\big)\cdot \delta q_{k}=0, 
\end{align}
where $t_1 = 0$ and $t_n = t_\textrm{f}$. 
The approximation can be equivalently represented by the forced discrete Euler-Lagrange (DEL) equation ($\delta q_1=0$ and $\delta q_n=0$ by definition):
\begin{align}
&\hspace{-0.3cm}D_1L_{\textrm{d}}(q_k,q_{k+1}) \hspace{-0.02cm}+\hspace{-0.02cm} D_2L_{\textrm{d}}(q_{k-1},q_{k}) +\nonumber\\&\hspace{-0.3cm}\ F^{-}_{\textrm{d}}(q_k,q_{k+1},u_k)\hspace{-0.02cm}+\hspace{-0.02cm} F^{+}_{\textrm{d}}(q_{k-1},q_{k},u_{k})=0. 
\end{align}
The DEL equation gives an implicit two-step mapping from two consecutive discrete configurations to the next $(q_{k-1},q_k)\rightarrow(q_{k+1})$.
Through the introduction of a new quantity, $p_k$, the forced DEL equations can be given in a \emph{position-momentum} form as
\begin{align}
&p_k + D_1L_{\textrm{d}}(q_k,q_{k+1})  + F^{-}_{\textrm{d}}(q_k,q_{k+1},u_k) = 0, \label{pm:form1}\\
&p_k =  D_2L_{\textrm{d}}(q_{k-1},q_{k}) + F^{+}_{\textrm{d}}(q_{k-1},q_{k},u_{k})\label{pm:form2}
\end{align}
and, therefore, an implicit one-step mapping is obtained $(q_{k},p_k)\rightarrow(q_{k+1})$\cite{West:thesis}.
Note that $p_k$ does not have any dependence on $q_{k+1}$.
Finally, by defining the integration equation 
\begin{align}
f(q_{k+1})=p_k\hspace{-0.08cm}+\hspace{-0.08cm}D_1L_{\textrm{d}}(q_k,q_{k+1})\hspace{-0.08cm}+\hspace{-0.08cm}F^{-}_{\textrm{d}}(q_k,q_{k+1},u_k),
\end{align}
and its derivative 
\begin{align}
Df(q_{k+1}) = D_2D_1L_{\textrm{d}} + D_2F^{-}_{\textrm{d}},
\end{align}
a simple root finding algorithm, outlined in Algorithm \ref{algo:1}, can be used to numerically solve the implicit relationship. 
Therefore, given $q_0$, $q_1$, and  the system's input, $u(t)$, the configuration can be propagated indefinitely. 
Furthermore, the variational integrator can be reformulated in the case where initial conditions are given in the form $(q(t_0),\dot{q}(t_0))$. 

\begin{algorithm}[H]
\caption{The Newton--Raphson Method}
\label{algo:1}
\begin{algorithmic}
\While {$|f(q_{k+1})|>\epsilon_{\textrm{tol}}$}
    \State{$q_{k+1} \leftarrow q_{k+1}  -Df^{-1}(q_{k+1})\cdot f(q_{k+1})$}
\EndWhile
\end{algorithmic}
\end{algorithm}

As discussed in \cite{VI:Murphey:02} and \cite{johnson2012trajectory} a first-order linearization of the discrete system dynamics can be obtained from the derived implicit one-step mapping (\ref{pm:form1})-(\ref{pm:form2}). 
The linearization is given in the following form
\begin{align}
\left[ \begin{array}{c} \delta q_{k+1}\\ \delta p_{k+1}\end{array} \right] &\hspace{-0.08cm}=\hspace{-0.08cm} \left[ \begin{array}{cc} \frac{\partial q_{k+1}}{\partial q_{k}}&\frac{\partial q_{k+1}}{\partial p_{k}}\\ \frac{\partial p_{k+1}}{\partial q_{k}}& \frac{\partial p_{k+1}}{\partial p_{k}}\end{array} \right]\left[ \begin{array}{c} \delta q_{k}\\ \delta p_{k}\end{array} \right]  \hspace{-0.08cm}+\hspace{-0.08cm} \left[ \begin{array}{c} \frac{\partial q_{k+1}}{\partial u_{k}}\\ \frac{\partial p_{k+1}}{\partial u_{k}}\end{array} \right] \delta u_k. \nonumber
\end{align}
The derivative $\frac{\partial q_{k+1}}{\partial q_k}$ is found by implicitly differentiating equation (\ref{pm:form1}):
\begin{align}
&\textstyle\frac{\partial}{\partial q_k}[p_k +  D_1L_{\textrm{d}}(q_k,q_{k+1})  + F^{-}_{\textrm{d}}(q_k,q_{k+1},u_k)=0], \nonumber \\ 
&\textstyle\frac{\partial q_{k+1}}{\partial q_k}= -M_{k+1}^{-1}\big( D_1D_1L_{\textrm{d}}(q_k,q_{k+1}) \nonumber\\&\quad\quad\quad\quad+ D_1F^{-}_{\textrm{d}}(q_k,q_{k+1},u_k)\big),  \label{SVI:partial1}
\end{align}
where $M_{k+1} =  D_2D_1L_{k+1} + D_2F^-_{k+1}$.
Quantities $\frac{\partial q_{k+1}}{\partial p_k}$ and $\frac{\partial q_{k+1}}{\partial u_k}$ can be found similarly. 
The remaining derivatives are found by explicitly differentiating equation (\ref{pm:form2}).
The same linearization procedure is applicable to constrained and stochastic systems \cite{VI:Murphey:02, johnson2012trajectory, IMAROB}. 

Recently, state estimation algorithms and numerical optimal control solvers have been shown to achieve better numerical performance when variational integrators, rather than one-step Euler methods, are used to linearize and propagate system dynamics \cite{DLT:DISS,IMAROB}.
However, traditional variational integrators can only achieve second-order accuracy and, therefore, further improvements cannot be easily achieved without reducing the time step--something that might be constrained by the embedded system or by computational resources. 
Section \ref{sec:backwards} introduces a method to increase the accuracy of variational integrators without increasing their complexity. 
Furthermore, numerical experiments in Section \ref{numex} demonstrate that a particular numerical optimal control solver performs better when fourth-order accurate variational integrators are utilized.

Typically, the left and right discrete forces are functions of the future input, $u_{k+1}$.
Specifically, the discrete forces are defined as
\begin{align}
&F^{\pm}_{\textrm{d}}(q_k,q_{k+1},u_{k})= \nonumber\\
&\quad\quad\quad\frac{h}{2}F(\frac{q_{k+1}+q_k}{2},\frac{q_{k+1}-q_{k}}{h}, \frac{u_k+u_{k+1}}{2}).
\end{align}
If the input signal is continuous this representation is required for the newly proposed surrogate variational integrator to obtain fourth-order convergence.    
However, if the discrete forces depend on the $u_{k+1}$ the first-order linearization given above becomes non-causal.   
Nevertheless, since the DDP algorithm produces a piecewise constant input signal the quantities $\frac{u_k+u_{k+1}}{2}$ and $u_k$ are equivalent since $u(T_k) = u_k, \ T_k \in [t_0+kh, \ t_0+(k+1)h )$.
That is, the average input over the time span of $T_k$ is $u_k$. 
Therefore, a causal linearization and a proper surrogate variational integrator can be defined in this setting.

\subsection{Overview of DDP}
The differential dynamic programming (DDP) algorithm iteratively solves a nonlinear optimal control problem using first- and second-order approximations of the considered dynamical system and the cost-to-go function, respectively \cite{DDPbook,DDPconstaint001,DDPminimax,SDDP}.
Specifically, the DDP algorithm is used to find the (local) minimum of a cost given by
\begin{align}\label{cost}
v(x,u,t) = h(x(t_\textrm{f}) + \int^{t_\textrm{f}}_{t_0} l(x(\tau), u(\tau), \tau) ~ \textrm{d}\tau,
\end{align}
subject to dynamics of the form
\begin{align}\label{dyna}
\dot x(t) = F(x(t),u(t)),
\end{align}
where $h(x(t_\textrm{f}))$ is the terminal cost, $l(x(t), u(t), t)$ is the running cost, $x$ is the state vector, and $u$ is the control input.
Being discrete in nature the DDP algorithm computes an optimal sequence of discrete inputs, $U=\{u_1,u_2,\dots\}$, to minimize the given cost such that the continuous input is then defined as $u(T_k) = u_k, \ T_k \in [t_0+kh, \ t_0+(k+1)h )$ where $h$ is the discretization time step.
Furthermore, the continuous system dynamics are also approximated by a user-selected numerical integrator as a sequence of state vectors such that $x_1 = x(t_0)$, $x_2 = x(t_0 + h)$, \dots, $x_N = x(t_\textrm{f})$. 
Iteratively the DDP algorithm updates the input, $U$, with a computed optimal control deviation, $U^\star$, as
\begin{align}\label{back:controlinput}
U \leftarrow U + \gamma\delta U^\star,
\end{align}
where $\gamma$ is a constant or generated with an automated process such as an Armijo line search \cite{Armi}. 
The optimization is then repeated using the newly updated input. 
Though the DDP algorithm has proven convergence guarantees several iterations can be required to obtained a solution sufficiently close to the optimal solution \cite{ShoemakerConv}. 

An abridged derivation of the DDP algorithm is given here and further details can be found in \cite{DDPbook,ShoemakerConv,SDDP}. 
To begin, given a nominal sequence of discrete inputs $\bar{U} = \{\bar{u}_1, \dots, \bar{u}_{n-1}\}$  a first-order linearization around the associated nominal state trajectory $\bar{X} =\{\bar{x}_1, \dots, \bar{x}_n\}$ is computed:
\begin{align} \label{equ:lin}
\delta x_{k+1}&=A_k \delta x_k + B_k \delta u_k,
\end{align}
where $A_k$ and $B_k$ are defined according to the utilized numerical integrator. 
Since variational integrators are considered in this paper the linearization defined in the previous section is used. 
Next, a second-order expansion of the optimal cost-to-go function
\begin{align} \label{eq:bellman}
V(\bar{x}_k,t_k) = \min_{u_k}\big[L(\bar{x}_k,\bar{u}_k,t_k) + V(\bar{x}_{k+1},t_{k+1})\big]
\end{align}
is given as \footnote{For ease of exposition, notation for derivatives is condensed to  $\nabla_{z}g = g_z$ and $\nabla_{xz}g = g_{xz}$. Furthermore, the condensed notation $Q_{x,k} =  Q_x(\bar{x}_k,\bar{u}_k)$ is used.} 
\begin{align}
&\min_{\delta u_k}\big[L(\bar{x}_k+\delta x_k,\bar{u}_k+\delta u_k,t_k) + V(\bar{x}_{k+1}+\delta x_{k+1},t_{k+1})\big] \nonumber \\
& \approx \min_{\delta u_k}\textstyle[Q_k + \delta u_k^\textrm{T} Q_{u,k} + \delta x_k^\textrm{T} Q_{x,k} + \frac{1}{2}\delta u_k^\textrm{T} Q_{uu,k} \delta u_k \nonumber\\
&\quad+ \frac{1}{2}\delta x_k^\textrm{T} Q_{xx,k} \delta x_k + \delta u_k^\textrm{T} Q_{ux,k} \delta x_k], \label{equ:expQ}
\end{align}
where 
\begin{align*}
L_k &= l_kh, \quad \quad \quad \quad \quad \quad \ \ Q_k  = \textstyle V_{k+1} + L_k,\\
Q_{x,k}  &= \textstyle L_{x,k} + A_k^\textrm{T}V_{x,k+1},\quad Q_{u,k}  =  \textstyle L_{u,k} + B_k^\textrm{T}V_{x,k+1} ,\\
Q_{xx,k}  &= \textstyle L_{xx,k} + A_k^\textrm{T}V_{xx,k+1} A_k,\\
Q_{uu,k}  &=  \textstyle L_{uu,k} + B_k^\textrm{T}V_{xx,k+1} B_k ,\\
Q_{xu,k}  &=  \textstyle L_{ux,k} + B_k^\textrm{T}V_{xx,k+1} A_k.
\end{align*}
Minimizing equation (\ref{equ:expQ}) yields the optimal control deviation $\delta U^\star = \{\delta u ^\star_1, \delta u ^\star_2, \dots\}$:
\begin{align} \label{back:optinput}
\delta u ^\star_k = -Q_{uu,k}^{-1}(Q_{u,k} + Q_{ux,k}\delta x^\star_k),
\end{align}
where $\delta X^\star = \{\delta x ^\star_1, \delta x ^\star_2, \dots\}$ is propagated as 
\begin{align}
\delta x^\star_k = A_k\delta x^\star_{k-1} + B_k\delta u^\star_{k-1}, \quad \delta x^\star_1 = 0.
\end{align}
The optimal control deviation contains a feedforward and a feedback component given by $Q_{uu,k}^{-1}Q_{u,k}$ and $Q_{uu,k}^{-1}Q_{ux,k}$, respectively.
A backward propagating second-order approximation of the value function is obtained when $\delta u ^\star$ is incorporated back into equation (\ref{equ:expQ}): 
\begin{align}
V_k &= V_{x, k+1} + L_k -  \frac{1}{2}Q_{u,k}Q_{uu,k}^{-1} Q_{u,k}, \label{prop:V}\\
V_{x,k} &= Q_{x,k} -  Q_{u,k}Q_{uu,k}^{-1} Q_{ux,k}^\textrm{T}, \label{prop:Vx}\\
V_{xx,k} &= Q_{xx,k} -  Q_{xu,k}Q_{uu,k}^{-1} Q_{ux,k}^\textrm{T},\label{prop:Vxx}
\end{align}
where the initial conditions are derived from the terminal cost, $h(\cdot)$, as $V(\bar{x}_{N}) = h(\bar{x}_N),V_x(\bar{x}_N) = h_x(\bar{x}_N)$, and $V_{xx}(\bar{x}_N) = h_{xx}(\bar{x}_N)$.
As stated earlier, the derived optimal control deviation, $\delta U ^\star$, is used to update the nominal input (\ref{back:controlinput}).
The optimization process can then be repeated to produce successively updated control inputs.
Termination of the optimization can occur when an update results in negligible change to cost or after a predetermined number of iterations.
The DDP algorithm is outlined in Algorithm \ref{back:algo:1}.

It should be expected that the accuracy of the utilized numerical integrator has a significant effect on the performance of the DDP algorithm. 
Note that both the forward propagation (discretized trajectory) and the backward propagation (value function approximation) rely on the selected numerical integrator. 
Therefore, effective optimized inputs cannot be obtained without a correct representation of the evolution of the dynamical system.

\begin{algorithm}[]
\caption{DDP with an Armijo Line Search}
\label{back:algo:1}
\begin{algorithmic}
\Require \\Initial discrete control input $u(t)$, parameters $\alpha, \beta, \epsilon$ \\system dynamics (\ref{dyna}), and cost function (\ref{cost}) 
\While {Cost updates results in more than $\epsilon$ in difference}
    \State Propagate the discretized trajectory
    \State Linearize the value function and system dynamics
    \State Back-propagate equations (\ref{prop:V})-(\ref{prop:Vxx})
    \State Compute $\delta U^\star$ and $\delta X^\star$ 
		\While {$\textrm{Cost}_\textrm{p} > \textrm{Cost} +\alpha\beta(v_x\delta X^\star + v_u\delta U^\star$)}
                      \State Find the proposed input $u_\textrm{p} \leftarrow u + \beta^j\delta U^\star$
		\State Propagate trajectory, $x_\textrm{p}$
                       \State Find proposed cost $\textrm{Cost}_\textrm{p}\leftarrow v(x_\textrm{p},u_\textrm{p},t)$
		\State Update $j \leftarrow j +1$ 
                      \EndWhile
 \State Update: $u \leftarrow u_\textrm{p}, x \leftarrow x_\textrm{p}, \textrm{Cost} \leftarrow \textrm{Cost}_\textrm{p}$
\EndWhile
\end{algorithmic}
\end{algorithm}

\section{Backward Error Analysis} \label{sec:backwards}
Backward error analysis quantifies the distortion induced to a differential equation by a particular numerical integrator. 
The methodology derives a \emph{modified differential equation} that exactly captures the propagated trajectory produced by the integrator. 
That is, by attempting to propagate a differential equation the error induced by a numerical method instead produces an error free sampled solution of the modified differential equation. 
Therefore, by comparing the differential equation and the modified differential equation the accuracy of a particular numerical integrator can be assessed directly in terms of the underlying state equations.
A rigorous treatment of backward error analysis can be found in References \cite{reich1999backward,hairer1994backward,WARMING1974159,griffiths1986scope}.  

To begin a brief summary of the methodology, suppose that a numerical integrator is used to approximate a continuous trajectory produced by an ordinary differential equation given by
 \begin{align}
\dot{x}(t) = f(x(t)), \quad x(0) = x_0, 
\end{align}
such that the discrete propagation $x_\textrm{d}=\{x_0,x_1,\dots,x_{n}\}$  is generated as 
 \begin{align}
x_{k+1} = \Psi(x_k),
\end{align}
where $x_k \approx x(k h)$ and $h$ is the discretization time step. 
Furthermore, consider a \emph{modified differential equation} of the form
 \begin{align}
f_{\textrm{mod}}(x(t)) &= f(x(t)) + h f_2(x(t))+ h^2 f_3(x(t)) \dots, \label{mod:diffeq}\\
\dot{\tilde{x}}(t) &= f_{\textrm{mod}}(\tilde{x}(t)), \quad \tilde{x}(0) = x_0, 
\end{align}
such that $x_k = \tilde{x}(k h)$. 
That is, the evolution of the modified differential equation is exactly captured by the numerical integrator. 
Under the assumption that the numerical method is expandable as
 \begin{align}
\Psi(x) = x + h f(x) + h^2\psi_2(x) + h^3\psi_3(x) + \dots,
\end{align} 
the modified differential equation is defined by a recursive relationship:
 \begin{align}
f_2(x) &= \psi_2(x) -\frac{1}{2!}f_{x}f, \\
f_3(x) &= \psi_3(x)-\frac{1}{3!}(f_{xx}\circ(f,f) + f_{x}f_{x}f) \nonumber\\&\quad-\frac{1}{2!}(f_{x}f_2 + f_{2,x}f), \\
f_4(x) &= \dots .
\end{align}
By examining the difference between $f(x(t))$ and $f_{\textrm{mod}}(x(t))$ the systematic error induced by the numerical integrator can be assessed.
Typically, the modified equation is a function of original  differential equation and the particular numerical integrator. 
Furthermore, the order of the numerical integrator determines if any terms of the modified differential equation are zero. 
For example, for a fourth-order integrator it is expected that $f_2(x)=0$ and $f_3(x)=0$. 

\subsection{Modified Lagrangians} \label{sec:backward:VI}
When analyzing variational integrators under the framework of backward error analysis \emph{modified Lagrangians} or \emph{Hamiltonians} are derived. 
That is, the modification caused by the variational integrator is characterized by another Hamiltonian system. 
A ``near-by" Hamiltonian system is exactly captured  by the variational integrator and, therefore, desirable geometric properties are maintained. 
This property of the numerical scheme offers insight to its energy and structure preserving abilities\cite{VI:Marsdem:01}. 

Reference \cite{vermeeren2015modified} presents a methodology for deriving a \emph{modified Lagrangian} when a midpoint variational integrator is used. 
The methodology supposes that there exist a modified Lagrangian, $L_\textrm{m}$, such that the discrete $L_\textrm{d}(q_k,q_{k+1})$, given by equation (\ref{dis:Lag}), equals the action integral of the modified system:
\begin{align}
 hL(\frac{q_{k+1}+q_k}{2},\frac{q_{k+1}-q_{k}}{h})=\int^{t_{k+1}}_{t_k}L_\textrm{m}(q,\dot{q}) \ \textrm{d}\tau.
\end{align}
After extensive analysis, the modified Lagrangian was found to be
\begin{align}\label{equ:modL1}
L_\textrm{m} = L &- \frac{h^2}{24}\big(-2\frac{\partial L}{\partial q}\ddot{q} + \dot{q}^\textrm{T}\frac{\partial^2 L}{\partial q\partial q}\dot{q}\nonumber\\&+ \ddot{q}^\textrm{T}\frac{\partial^2 L}{\partial \dot q\partial \dot q}\ddot{q}+ 2\ddot{q}^\textrm{T}\frac{\partial^2 L}{\partial q\partial \dot q}\dot{q}\big) + o(h^4).  
\end{align}
The equations of motion of the modified and original  systems can be related as
\begin{align}\label{equ:modsys}
&\ddot{q}(t) = \Big(\frac{\partial^2 L_\textrm{m}}{\partial \dot q\partial \dot q}\Big)^{-1}\Big(\frac{\partial L_\textrm{m}}{\partial q} - \frac{\partial^2 L_\textrm{m}}{\partial q\partial \dot q}\dot q\Big)\nonumber\\&\quad\quad=\Big(\frac{\partial^2 L}{\partial \dot q\partial \dot q}\Big)^{-1}\Big(\frac{\partial L}{\partial q} - \frac{\partial^2 L}{\partial q\partial \dot q}\dot q\Big)+ o(h^2).
\end{align}
since the respective Lagrangians differ in terms of order $h^2$.
Furthermore, the equations of motion (\ref{equ:modsys}) can be incorporated into the expression of the modified Lagrangian to eliminate any dependence on $\ddot{q}(t)$:
\begin{align} \label{equ:modL}
L_\textrm{m} &=  {L} -\frac{h^2}{24}\Big(\dot{q}^T\big(\frac{\partial^2 {L}}{\partial q\partial q} - \frac{\partial^2 {L}}{\partial\dot q\partial q}^\textrm{T}\frac{\partial^2 {L}}{\partial\dot q\partial\dot q}^{-1}\frac{\partial^2 {L}}{\partial\dot q\partial q}\big)\dot{q} \nonumber\\
&- \frac{\partial {L}}{q}^\textrm{T}\frac{\partial^2 {L}}{\partial\dot q\partial\dot q}^\textrm{-1}\frac{\partial {L}}{q}
+ 2\dot{q}^\textrm{T}\frac{\partial^2 {L}}{\partial\dot q\partial q}^\textrm{T}\frac{\partial^2{L}}{\partial\dot q\partial\dot q}^{-1}\frac{\partial {L}}{\partial q}\Big) + o(h^4).\nonumber
\end{align}

\subsection{Surrogate Lagrangians}
In the previous section we have quantified the error induced by a variational integrator.
It is only natural to look for a manner in which to mitigate or eliminate the quantified error. 
Several numerical integrators have been developed for specific classes of dynamical systems that eliminate integration error (up to some order)
 \cite{mushtaq2014higher,Mushtaq20141461,chartier2007numerical,hairer2006preprocessed,StochModEqu,chartier2007modified,kozlov2008high}.
However, the surrogate variational integrator outlined below is unique in that no assumption on the structure of the system's Lagrangian is needed.

The proposed approach uses the same integration scheme as traditional variational integrators, but replaces the Lagrangian of the system with a surrogate.
Through the introduction of the surrogate Lagrangian the induced error of order $h^2$ can be eliminated and, therefore, an increase in accuracy is achieved. 
The \emph{surrogate} Lagrangian, $\hat{L}$, is defined as
\begin{align} 
\hat{L} =  {L} &+ \frac{h^2}{24}(-2\frac{\partial L}{\partial q}\ddot{q} + \dot{q}^\textrm{T}\frac{\partial^2 L}{\partial q\partial q}\dot{q}\nonumber\\ 
&+ \ddot{q}^\textrm{T}\frac{\partial^2 L}{\partial \dot q\partial \dot q}\ddot{q}+ 2\ddot{q}^\textrm{T}\frac{\partial^2 L}{\partial q\partial \dot q}\dot{q}).
\end{align}
The \emph{ modified surrogate} Lagrangian, $\hat{L}_\textrm{m}$,  is derived from equations (\ref{equ:modL1}) and (\ref{equ:modL}) as 
\begin{align} 
\hat{L}_\textrm{m} =  \hat{L} &- \frac{h^2}{24}(-2\frac{\partial \hat{L}}{\partial q}\ddot{q} + \dot{q}^\textrm{T}\frac{\partial^2 \hat{L}}{\partial q\partial q}\dot{q}\nonumber\\&+ \ddot{q}^\textrm{T}\frac{\partial^2 \hat{L}}{\partial \dot q\partial \dot q}\ddot{q}+ 2\ddot{q}^\textrm{T}\frac{\partial^2 \hat{L}}{\partial q\partial \dot q}\dot{q}).
\end{align}
Therefore, 
\begin{align}
\hat{L}_\textrm{m}  &=  {L} + o(h^4), \label{equ:mod2}
\end{align}
since
\begin{align} 
&-2\frac{\partial \hat{L}}{\partial q}\ddot{q} + \dot{q}^\textrm{T}\frac{\partial^2 \hat{L}}{\partial q\partial q}\dot{q}+ \ddot{q}^\textrm{T}\frac{\partial^2 \hat{L}}{\partial \dot q\partial \dot q}\ddot{q}+ 2\ddot{q}^\textrm{T}\frac{\partial^2 \hat{L}}{\partial q\partial \dot q}\dot{q}+ o(h^2) \nonumber\\&\quad= -2\frac{\partial L}{\partial q}\ddot{q} + \dot{q}^\textrm{T}\frac{\partial^2 L}{\partial q\partial q}\dot{q}+ \ddot{q}^\textrm{T}\frac{\partial^2 L}{\partial \dot q\partial \dot q}\ddot{q}+ 2\ddot{q}^\textrm{T}\frac{\partial^2 L}{\partial q\partial \dot q}\dot{q}.
\end{align}
Note that (\ref{equ:mod2}) implies that the difference in equations of motion of the modified surrogate and original systems is of order $h^4$. 
Therefore, the surrogate variational integrator, a second-order integrator, achieves fourth-order accuracy. 
No change in the central integration scheme was needed in order to obtained the increase in accuracy. 
Furthermore, since the central integration scheme remains the same a linearization of the discrete system dynamics is available. 
Therefore, the surrogate variational integrator can be used in algorithms requiring linearization of system dynamics (e.g. numerical optimal control solvers\cite{rao2009survey}, Kalman filters\cite{simon2006optimal}, etc.). 
Lastly, the methodology presented in Reference \cite{vermeeren2015modified} is applicable to forced and constrained systems.
Therefore, surrogate Lagrangians for constrained and forced systems can be derived, but due to space constraints further analysis is not given. 
However, the numerical examples given in Section \ref{numex} are of forced systems. 

\section{Numerical Examples} \label{numex}
In this section, DDP algorithms utilizing the nominal and surrogate variational integrators are used to numerically solve two optimal control problems. 
It is shown that the DDP algorithms utilizing the surrogate variational integrators are able to produce more optimized inputs and trajectories. 
Furthermore, the accuracy of the trajectories produced by the variational integrators are also compared. 
As expected, the surrogate variational integrator achieves fourth-order convergence in propagation error. 

Two metrics are used to judge the performance of the variational integrators: propagation error $\mathcal{L}^2$-norm and assessed cost. 
The  \emph{propagation error $\mathcal{L}^2$-norm} of a variational integrator is computed as
\begin{align}\label{prop:error}
e_{L^2}  &= \Big(h\sum^{t_n}_{t_0}(q_k - q_\textrm{bm}(t_k))^2 \Big)^{\frac{1}{2}}
\end{align}
where $q_\textrm{bm}(t)$ is the trajectory obtained from a high-fidelity benchmark integrator (analytic or sufficiently small discretization time step). 
The \emph{assessed cost} is the cost incurred by the system when the optimized piecewise constant  inputs are used in a continuous-time setting. 
For the presented experiments, the continuous response of the system is approximated by propagating the system with a high-fidelity benchmark integrator. 
Therefore, the assessed cost is able to appraise the effect of optimizing the continuous problem in a discrete setting.  
In particular, the assessed cost will change as the discretization time step used by the DDP algorithm is varied. 

\subsection{Harmonic Oscillator}

In this example a harmonic oscillator is considered with mass $M$ and spring constant $K$.
Its Lagrangian is given as 
\begin{eqnarray}
{L}&=  \frac{1}{2}M\dot q^2 - \frac{1}{2}K q^2,
\end{eqnarray}
and its surrogate Lagrangian is given as
\begin{eqnarray}
\hat{L}&=  \frac{1}{2}(M-\frac{h^2}{12}K)\dot q^2 - \frac{1}{2}(K+\frac{h^2}{12}K{M}^{-1}K)q^2,
\end{eqnarray}
where
 \[
M =
\begin{bmatrix}
    2& 0.1 & 0 &  0.3 \\
 0.1&3&0.1&0\\
0&0.1&4.1&0.3\\  
0.3&0&0.3&4
\end{bmatrix},
\]
 \[
K = 
\begin{bmatrix}
1&0.5&0&0.5\\ 
0.5&0.9&0.35&0\\ 
0&0.35&8.1&0.65\\  
0.5&0&0.65&2.1
\end{bmatrix}.
\]
The systems equations of motion are given as
\begin{eqnarray}
M\ddot{q}&=  -Kq+Bu(t), \quad q(t_0) = q_0, \quad \dot q(t_0) = \dot q_0,
\end{eqnarray}
where  \[
B = 
\begin{bmatrix}
1&0\\ 
0&0\\ 
0&1\\  
0&0
\end{bmatrix},
\] is the input matrix,
$u$ is the input, and
$q_0 = [5,0,0,0]^\textrm{T}$,
$\dot q_0 = [0,0,0,0]^\textrm{T}$ are the initial conditions.
Note that the system is under-actuated and is not damped. 
Therefore, an accurate prediction of future states is critical to determine an optimal control policy.
The cost to be minimized is
\begin{align}
v(x,u,t) = \int^{10}_{0} \frac{1}{20}u(\tau)^\textrm{T}u(\tau) + \frac{1}{2}q(\tau)^\textrm{T}q(\tau) ~ \textrm{d}\tau.
\end{align}
Recall that the DDP algorithm is discrete in nature and, therefore, a piecewise constant, not a continuous, input is optimized.  
The piecewise input is defined as $u(T_k) = u_k, \ T_k \in [t_0+kh, \ t_0+(k+1)h )$ where $h$ is the discretization time step.
Finally, since the considered system is linear the configuration can be propagated error free using an explicit integrator:
 \begin{align}\label{exact:prop}
x_{k+1} = e^{Ah}x_k + A^{-1}(e^{Ah}-\textrm{I})Bu_k,
\end{align}
where $A$ is the state transition matrix of the system in its state space formulation.
Furthermore, note that (\ref{exact:prop}) also provides a linearization of the discrete system dynamics. 
The numerical integrator given in (\ref{exact:prop}) is the \emph{benchmark integrator} for this example.
Of course, for linear systems the best numerical integrator to use in the DDP algorithm would be (\ref{exact:prop}). 
Therefore, the variational integrators can be fairly assessed by being compared to  (\ref{exact:prop}). 

\begin{figure}
    \centering
        	\includegraphics[width=0.9\columnwidth]{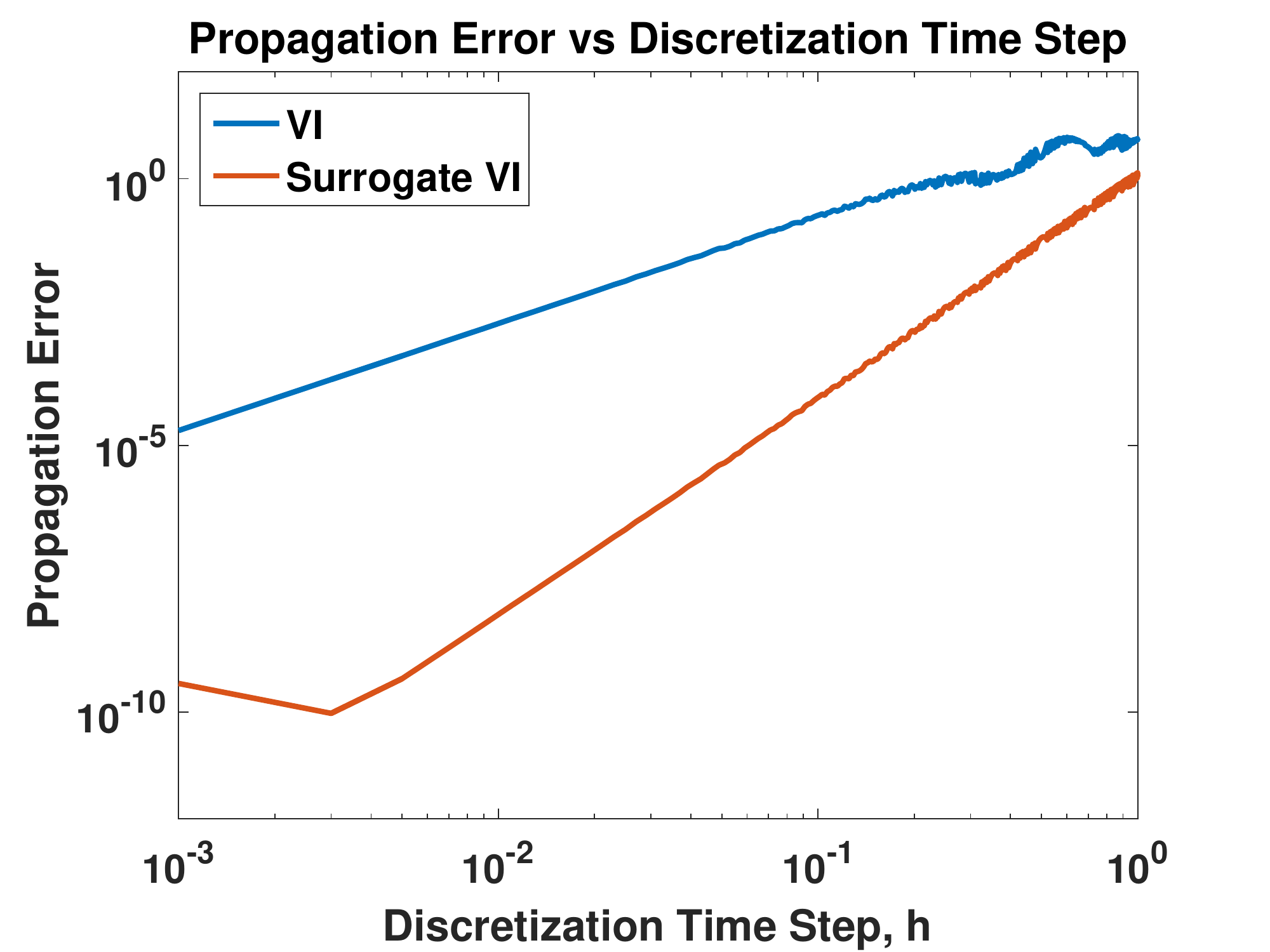}
    \caption{ The propagation error $L^2$-norm (benchmark trajectory obtained by (\ref{exact:prop})) as a function of the utilized discretization time step. 
The nominal and surrogate variational integrators display second- and fourth-order convergence, respectively. 
}
    \label{fig1:ex1}
\end{figure}

Figure \ref{fig1:ex1} shows the propagation error $\mathcal{L}^2$-norm as a function of the discretization time step when the harmonic oscillator was simulated for 250 seconds and all inputs were set to unity.  
As predicted by the analysis outlined in Section \ref{sec:backwards} the surrogate variational integrator obtained a fourth-order convergence of the propagation error. 

\begin{figure}
	\centering
        	\includegraphics[width=0.9\columnwidth]{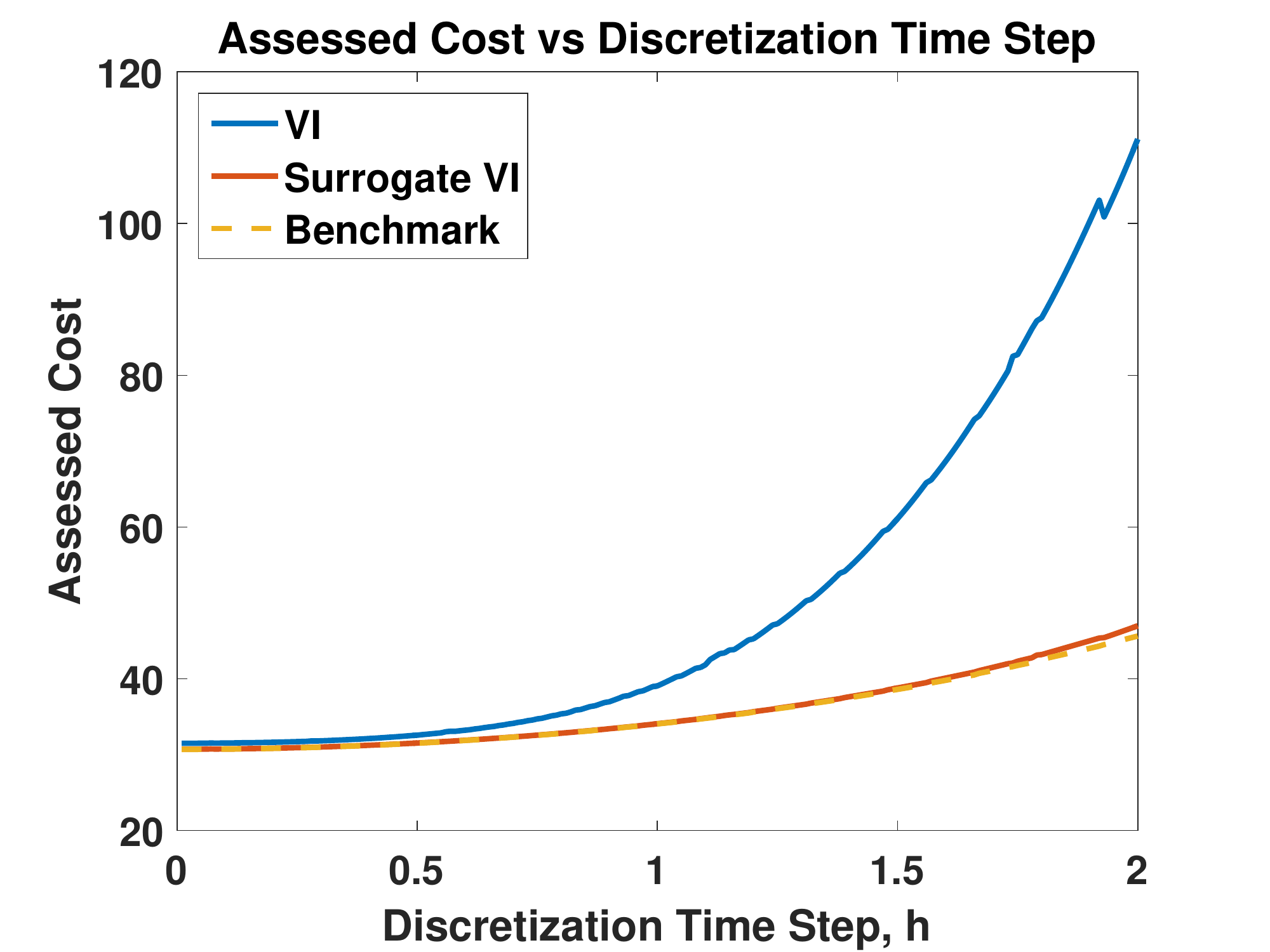}
    \caption{The assessed cost of the optimized control inputs and associated trajectories as a function of the utilized discretization time step.
The assessed cost grew significantly faster as a function of the discretization time step when the nominal variational integrator was used.
The surrogate variational and benchmark integrators (\ref{exact:prop}) resulted in very similar assessed costs. 
}
    \label{fig2:ex1}
\end{figure}

Figure \ref{fig2:ex1} shows the assessed cost of optimized control inputs as a function of the utilized discretization time step for the three numerical integrators considered.
The assessed cost grew significantly faster as a function of the discretization time step when the nominal variational integrator was used.
However, the surrogate variational and benchmark integrators produced very similar assessed costs.
Note that the assessed cost for all integrators increases as the discretization time step increases.  
This increase can be attributed to a decrease in controller bandwidth, an increase in integration errors, and numerical errors arising from approximations found in the derivation of the DDP algorithm. 
However, the increase seen for the benchmark integrator is not caused by errors in propagation or linearization.
Therefore, it can be concluded that the surrogate variational integrator is sufficiently accurate to eliminate most affects of integration error on the DDP algorithm.   

\subsection{Cart with Double Pendulum}
\begin{figure}
	\centering
        	\includegraphics[width=0.8\columnwidth]{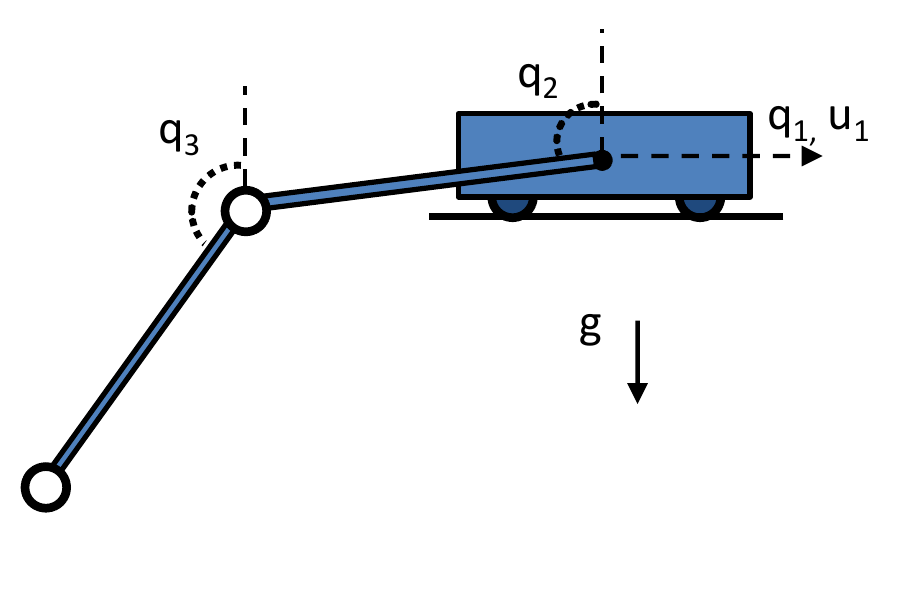}
    \caption{Diagram of the considered cart with a double pendulum system.
The system consists of a single input and three generalized coordinates where $q_1$ is the horizontal position of the cart, and $q_2$ and $q_3$ are the rotation angles with respect to the inertial frame of the first and second pendulum, respectively.
}
    \label{fig2:dia}
\end{figure}

A cart with a double pendulum attached to its center of mass, as depicted in Figure \ref{fig2:dia}, is considered for this example (see \cite{Graichen200763} for a full description of this dynamical system). 
The mass of the cart and both pendulums are 1 kilogram and the lengths of the first and second pendulum are 4 and 6 meters, respectively. 
The mass of the pendulums are assumed to be concentrated at their ends.
The system was subjected to a gravitational field ($\textrm{9.81 m/s}^2$) and no damping forces were modeled. 
There is a single input that directly influences the acceleration of the cart such that the systems equations of motion are given as
\begin{align}
\frac{\partial^2L(q,\dot{q})}{\partial^2\dot{q}}\ddot{q}&=  \frac{\partial L(q,\dot{q})}{\partial{q}}-\frac{\partial^2L(q,\dot{q})}{\partial{q}\partial{\dot q}}\dot q+Bu(t), \nonumber\\ &\quad q(t_0) = q_0, \quad \dot q(t_0) = \dot q_0
\end{align}
where  $L$ is the system Lagrangian, $B =[1,0,0]^\textrm{T}$, $q(t_0) = [0,\frac{9\pi}{8},\frac{5\pi}{6}]^\textrm{T}$, $\dot{q}(t_0) = [0,0,0]^\textrm{T}$, $q_1$ is the horizontal position of the cart in meters, and $q_2$ and $q_3$ are the rotation angles in radians with respect to the inertial frame of the first and second pendulum, respectively.
If $q = 0$ and $\dot{q}=0$ the double pendulum is inverted and in an unstable equilibrium. 
The cost to be minimized is
\begin{align}
v(x,u,t) &= \int^{20}_{0} \frac{1}{200}u(\tau)^\textrm{T}u(\tau) +50q_1(\tau)^2 \nonumber\\&\quad+ 50\pi \big(q_2(\tau)-\pi)^2 + (q_3(\tau)-\pi)^2\big) ~ \textrm{d}\tau.
\end{align}
Therefore, the goal is to drive the system to its stable equilibrium. 
Since the system is nonlinear the benchmark integrator used in the previous example is not available. 
Therefore, the benchmark integrator is the surrogate variational integrator when $h=1\times10^{-3}$. 

Unlike the previous example, the surrogate Lagrangian cannot be parameterized by a change in physical parameters. 
Several more polynomials appear in the surrogate Lagrangian than in the nominal Lagrangian. 
Therefore, evaluation of the discrete surrogate Lagrangian and its derivatives is computationally more expensive than in the nominal case.
However, the complexity of the central integration scheme remains the same. 
The added expense of evaluating surrogate Lagrangians can be mitigated through parallelization or code optimization.  

\begin{figure}
    \centering
        	\includegraphics[width=0.9\columnwidth]{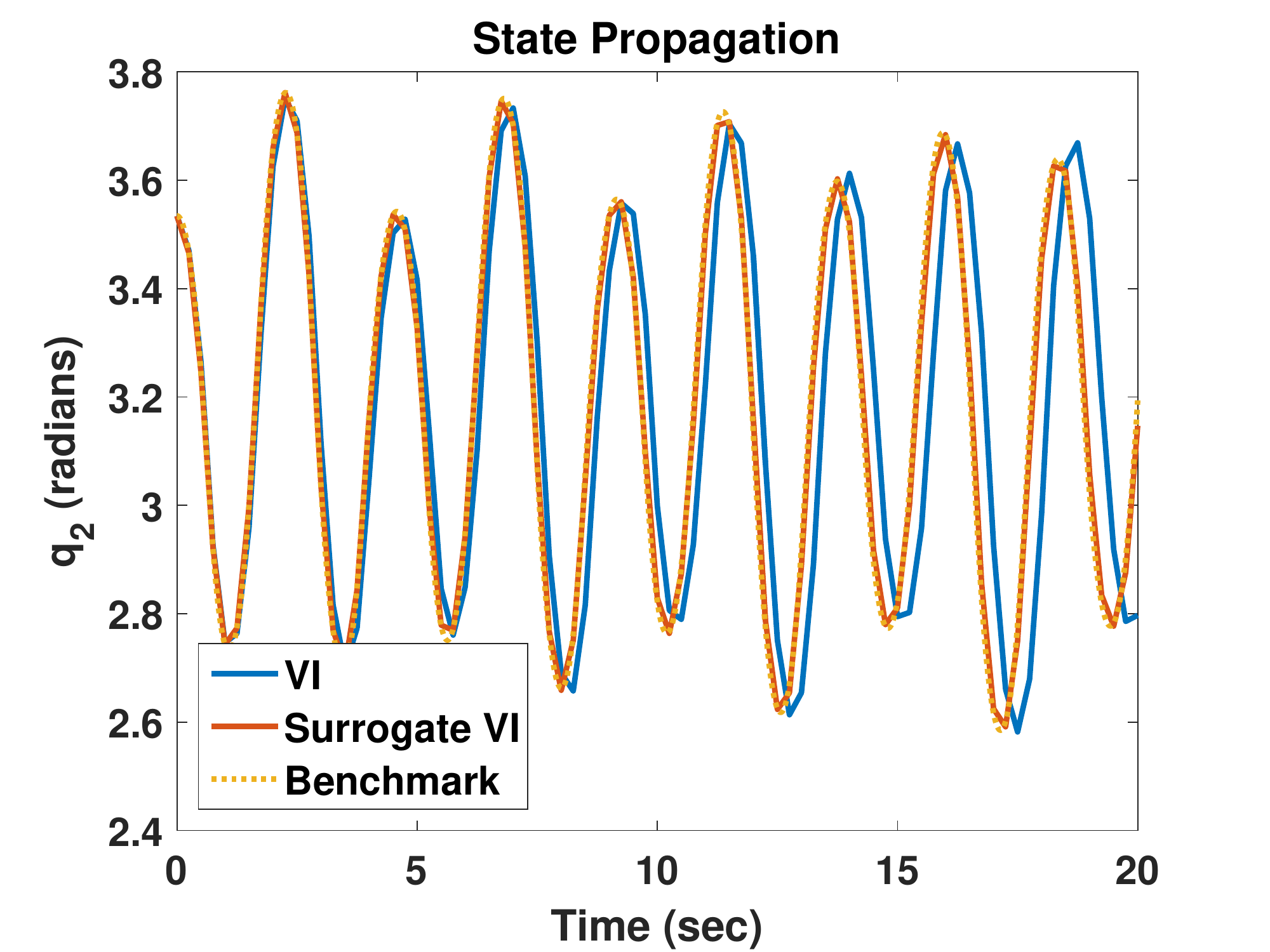}
    \caption{The propagated $q_2$ trajectory obtained by the nominal and surrogate variational integrator when $h=0.25$ and all inputs were set to unity. 
The benchmark trajectory is obtained by the surrogate variational integrator when $h=1\times10^{-3}$.
While the nominal variational integrator produces significant propagation errors, the surrogate variational integrator is able to accurately predict the evolution of the system.  
}
    \label{fig1:ex2}
\end{figure}

Figure \ref{fig1:ex2} shows the propagated $q_2$ trajectory obtained by the nominal and surrogate variational integrator when $h=0.25$ and all inputs were set to unity. 
While the nominal variational integrator produces significant propagation errors, the surrogate variational integrator is able to accurately predict the evolution of the system.  
The difference of the propagated $q_2$ trajectories and the benchmark trajectory at $t=20$ was $0.3967$ and $0.0470$ when the nominal and surrogate variational integrator were used, respectively. 
Although not shown due to space constraints, the surrogate variational integrator obtained fourth-order convergence of the propagation error. 

\begin{figure}
	\centering
        	\includegraphics[width=0.9\columnwidth]{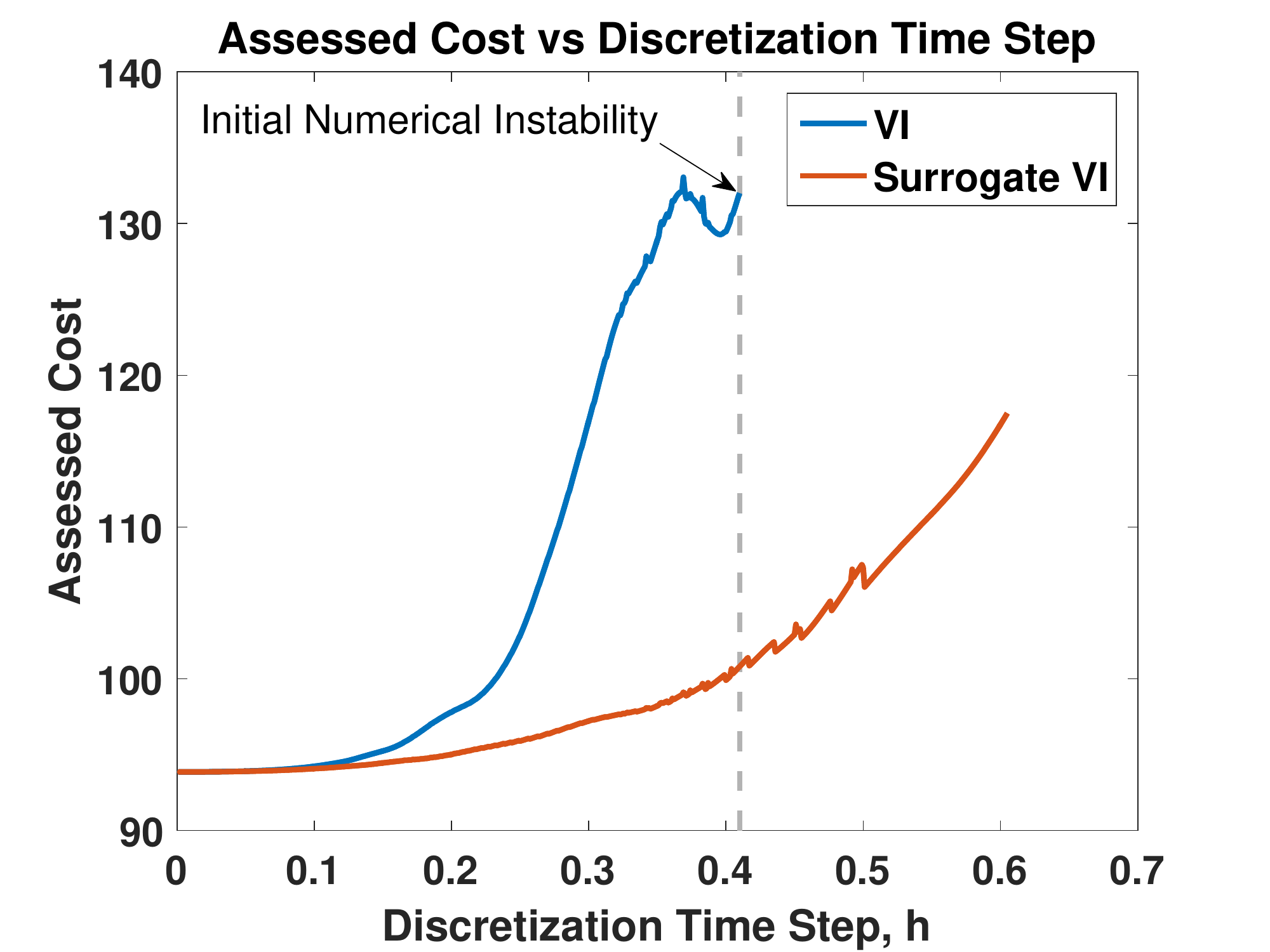}
    \caption{The assessed cost of the optimized control inputs and associated trajectories as a function of the utilized discretization time step.
The assessed cost grew significantly faster as a function of the discretization time step when the nominal variational integrator was used.
On the other hand, the assessed cost was relatively constant when the surrogate variational integrator was used. 
While computing the assessed cost, the solutions produced by the nominal and surrogate variational integrators initially became unstable when $h=0.410$ and $h=0.605$, respectively. 
}
    \label{fig2:ex2}
\end{figure}

Figure \ref{fig2:ex2} shows the assessed cost of optimized control inputs as a function of the utilized discretization time step. 
The assessed cost grew significantly faster as a function of the discretization time step when the nominal variational integrator was used.
Qualitatively, Figures \ref{fig2:ex1} and \ref{fig2:ex2} show similar trends.
As before, the difference in numerical integrators result in a significant difference in the performance of the overall DDP algorithm.
While computing the assessed cost, the solutions produced by the nominal and surrogate variational integrators initially became unstable when $h=0.410$ and $h=0.605$ was used in the DDP algorithm, respectively.  
At these discretization time steps the assessed cost experienced an abrupt change of multiple orders of magnitude.
This abrupt change suggests that the DDP algorithm failed to converged to a useful solution.

\section{Conclusion}
This paper introduced the concept of a surrogate variational integrator and demonstrated the benefits of using the proposed integrator in the DDP algorithm. 
These benefits are not expected to be exclusive to the DDP algorithm.
Any algorithm or process that utilizes propagated system trajectories or linearizes system dynamics will undoubtedly benefit. 
The presented work motivates further investigation into the numerical properties of algorithms used in control and state estimation. 
Furthermore, backward error analysis may provide the necessary tools to analytically determine an algorithm's dependence on computational resources.
Therefore, the design process of control and state estimation algorithms can explicitly take into account, for example, constraints on memory allocation and execution time. 
While the complexity of the central variational integration scheme is not increased by the introduction of surrogate Lagrangians, evaluation of discrete surrogate Lagrangian and its derivatives can be computationally expensive. 
However, the computational cost of these evaluations can be mitigated through parallelization or code optimization.  
\bibliographystyle{IEEEtran}
\bibliography{references1,references2}

\begin{thebibliography}{10}
\providecommand{\url}[1]{#1}
\csname url@samestyle\endcsname
\providecommand{\newblock}{\relax}
\providecommand{\bibinfo}[2]{#2}
\providecommand{\BIBentrySTDinterwordspacing}{\spaceskip=0pt\relax}
\providecommand{\BIBentryALTinterwordstretchfactor}{4}
\providecommand{\BIBentryALTinterwordspacing}{\spaceskip=\fontdimen2\font plus
\BIBentryALTinterwordstretchfactor\fontdimen3\font minus
  \fontdimen4\font\relax}
\providecommand{\BIBforeignlanguage}[2]{{%
\expandafter\ifx\csname l@#1\endcsname\relax
\typeout{** WARNING: IEEEtran.bst: No hyphenation pattern has been}%
\typeout{** loaded for the language `#1'. Using the pattern for}%
\typeout{** the default language instead.}%
\else
\language=\csname l@#1\endcsname
\fi
#2}}
\providecommand{\BIBdecl}{\relax}
\BIBdecl

\bibitem{rao2009survey}
A.~V. Rao, ``A survey of numerical methods for optimal control,''
  \emph{Advances in the Astronautical Sciences}, vol. 135, no.~1, pp. 497--528,
  2009.

\bibitem{DDPbook}
D.~H. Jacobson and D.~Q. Mayne, \emph{Differential Dynamic Programming}.\hskip
  1em plus 0.5em minus 0.4em\relax Elsevier, 1970.

\bibitem{DDPconstaint001}
G.~Lantoine and R.~P. Russell, ``A hybrid differential dynamic programming
  algorithm for robust low-thrust optimization,'' in \emph{AAS/AIAA
  Astrodynamics Specialist Conference and Exhibit}, 2008.

\bibitem{DDPminimax}
J.~Morimoto and C.~G. Atkeson, ``Minimax differential dynamic programming: An
  application to robust biped walking,'' in \emph{International Conference on
  Intelligent Robots and Systems}, 2003, pp. 1927--1932.

\bibitem{SDDP}
E.~Theodorou, Y.~Tassa, and E.~Todorov, ``Stochastic differential dynamic
  programming,'' in \emph{American Control Conference}, 2010, pp. 1125--1132.

\bibitem{ilqr001}
E.~Todorov and W.~Li, ``A generalized iterative lqg method for locally-optimal
  feedback control of constrained nonlinear stochastic systems,'' in
  \emph{American Control Conference}, 2005, pp. 300--306.

\bibitem{DDPconstaint002}
S.~J. Yakowitz, ``The stagewise kuhn-tucker condition and differential dynamic
  programming,'' \emph{IEEE Transactions on Automatic Control}, vol.~31, no.~1,
  pp. 25--30, 1986.

\bibitem{ShoemakerConv}
L.~Z. Liao and C.~A. Shoemaker, ``Convergence in unconstrained discrete-time
  differential dynamic programming,'' \emph{IEEE Transactions on Automatic
  Control}, vol.~36, no.~6, pp. 692--706, 1991.

\bibitem{DLT:DISS}
G.~De~La~Torre, ``Autonomous suspended load operations via trajectory
  optimization and variational integrators,'' Ph.D. dissertation, Georgia
  Institute of Technology, 2015.

\bibitem{IMAROB}
G.~De~La~Torre and E.~Theodorou, ``Stochastic variational integrators for
  system propagation and linearization,'' in \emph{IMA Conference on
  Mathematics of Robotics}, 2015.

\bibitem{ober2011discrete}
S.~Ober-Bl{\"o}baum, O.~Junge, and J.~E. Marsden, ``Discrete mechanics and
  optimal control: an analysis,'' \emph{ESAIM: Control, Optimisation and
  Calculus of Variations}, vol.~17, no.~2, pp. 322--352, 2011.

\bibitem{VI:Marsdem:01}
J.~E. Marsden and M.~West, ``Discrete mechanics and variational integrators,''
  \emph{Acta Numerica}, vol.~10, pp. 357--514, 2001.

\bibitem{VI:Murphey:02}
E.~Johnson, J.~Schultz, and T.~D. Murphey, ``Structured linearization of
  discrete mechanical systems for analysis and optimal control,'' \emph{IEEE
  Transactions on Automation Science and Engineering}, vol.~12, no.~1, pp.
  140--152, 2015.

\bibitem{hussein2006discrete}
I.~I. Hussein, M.~Leok, A.~K. Sanyal, and A.~M. Bloch, ``A discrete variational
  integrator for optimal control problems on {SO(3)},'' in \emph{IEEE
  Conference on Decision and Control}, 2006, pp. 6636--6641.

\bibitem{johnson2012trajectory}
E.~Johnson, ``Trajectory optimization and regulation for constrained discrete
  mechanical systems,'' Ph.D. dissertation, Northwestern University, 2012.

\bibitem{reich1999backward}
S.~Reich, ``Backward error analysis for numerical integrators,'' \emph{SIAM
  Journal on Numerical Analysis}, vol.~36, no.~5, pp. 1549--1570, 1999.

\bibitem{hairer1994backward}
E.~Hairer, ``Backward analysis of numerical integrators and symplectic
  methods,'' \emph{Annals of Numerical Mathematics}, vol.~1, pp. 107--132,
  1994.

\bibitem{WARMING1974159}
R.~F. Warming and B.~J. Hyett, ``The modified equation approach to the
  stability and accuracy analysis of finite-difference methods,'' \emph{Journal
  of Computational Physics}, vol.~14, no.~2, pp. 159--179, 1974.

\bibitem{griffiths1986scope}
D.~F. Griffiths and J.~M. Sanz-Serna, ``On the scope of the method of modified
  equations,'' \emph{SIAM Journal on Scientific and Statistical Computing},
  vol.~7, no.~3, pp. 994--1008, 1986.

\bibitem{West:thesis}
M.~West, ``Variational integrators,'' Ph.D. dissertation, California Institute
  of Technology, 2004.

\bibitem{VI:Murphey:01}
E.~R. Johnson and T.~D. Murphey, ``Scalable variational integrators for
  constrained mechanical systems in generalized coordinates,'' \emph{IEEE
  Transactions on Robotics}, vol.~25, no.~6, pp. 1249--1261, 2009.

\bibitem{2002analytical}
A.~I. Lurie, \emph{Analytical Mechanics}.\hskip 1em plus 0.5em minus
  0.4em\relax Springer, 2002.

\bibitem{Armi}
C.~T. Kelley, \emph{Iterative methods for optimization}.\hskip 1em plus 0.5em
  minus 0.4em\relax SIAM, 1999.

\bibitem{vermeeren2015modified}
M.~Vermeeren, ``Modified equations for variational integrators,'' \emph{arXiv},
  2015.

\bibitem{mushtaq2014higher}
A.~Mushtaq, A.~Kv{\ae}rn{\o}, and K.~Olaussen, ``Higher-order geometric
  integrators for a class of {H}amiltonian systems,'' \emph{International
  Journal of Geometric Methods in Modern Physics}, vol.~11, no.~1, pp.
  1\,450\,009--1–1\,450\,009--20, 2014.

\bibitem{Mushtaq20141461}
A.~Mushtaq and K.~Olaussen, ``Automatic code generator for higher order
  integrators,'' \emph{Computer Physics Communications}, vol. 185, no.~5, pp.
  1461--1472, 2014.

\bibitem{chartier2007numerical}
P.~Chartier, E.~Hairer, and G.~Vilmart, ``Numerical integrators based on
  modified differential equations,'' \emph{Mathematics of Computation},
  vol.~76, no. 260, pp. 1941--1953, 2007.

\bibitem{hairer2006preprocessed}
E.~Hairer and G.~Vilmart, ``Preprocessed discrete {M}oser--{V}eselov algorithm
  for the full dynamics of a rigid body,'' \emph{Journal of Physics A:
  Mathematical and General}, vol.~39, no.~42, pp. 13\,225--13\,235, 2006.

\bibitem{StochModEqu}
A.~Abdulle, D.~Cohen, G.~Vilmart, and K.~C. Zygalakis, ``High weak order
  methods for stochastic differential equations based on modified equations,''
  \emph{SIAM Journal on Scientific Computing}, vol.~34, no.~3, pp. 1800--1823,
  2012.

\bibitem{chartier2007modified}
{Chartier, Philippe}, {Hairer, Ernst}, and {Vilmart, Gilles}, ``Modified
  differential equations,'' \emph{ESAIM: PROCEEDINGS}, vol.~21, pp. 16--20,
  2007.

\bibitem{kozlov2008high}
R.~Kozlov, ``High-order conservative discretizations for some cases of the
  rigid body motion,'' \emph{Physics Letters A}, vol. 373, no.~1, pp. 23--29,
  2008.

\bibitem{simon2006optimal}
D.~Simon, \emph{Optimal state estimation: Kalman, {H$_{infinity}$}, and
  nonlinear approaches}.\hskip 1em plus 0.5em minus 0.4em\relax John Wiley \&
  Sons, 2006.

\bibitem{Graichen200763}
K.~Graichen, M.~Treuer, and M.~Zeitz, ``Swing-up of the double pendulum on a
  cart by feedforward and feedback control with experimental validation,''
  \emph{Automatica}, vol.~43, no.~1, pp. 63--71, 2007.

\end{thebibliography}
\end{document}